\documentclass[11pt,leqno]{article}
\usepackage{amsthm,amsfonts,amssymb,amsmath,oldgerm}
\usepackage{epsfig}
\numberwithin{equation}{section}
\usepackage[thinlines]{easybmat}

\renewcommand\b{\beta}

\renewcommand\b{\beta}

\newcommand\br{\begin{remark}}
\newcommand\er{\end{remark}}
\newcommand\bp{\begin{pmatrix}}
\newcommand\ep{\end{pmatrix}}
\newcommand\be{\begin{equation}}
\newcommand\ee{\end{equation}}
\newcommand\ba{\begin{equation}\begin{aligned}}
\newcommand\ea{\end{aligned}\end{equation}}

\newcommand{\bap}{\begin{app}}
\newcommand{\eap}{\end{app}}
\newcommand{\begs}{\begin{exams}}
\newcommand{\eegs}{\end{exams}}
\newcommand{\beg}{\begin{example}}
\newcommand{\eeg}{\end{exaplem}}
\newcommand{\bpr}{\begin{proposition}}
\newcommand{\epr}{\end{proposition}}
\newcommand{\bt}{\begin{theorem}}
\newcommand{\et}{\end{theorem}}
\newcommand{\bc}{\begin{corollary}}
\newcommand{\ec}{\end{corollary}}
\newcommand{\bl}{\begin{lemma}}
\newcommand{\el}{\end{lemma}}
\newcommand{\bd}{\begin{definition}}
\newcommand{\ed}{\end{definition}}
\newcommand{\brs}{\begin{remarks}}
\newcommand{\ers}{\end{remarks}}

\newtheorem{theo}{Theorem}[section]

\newtheorem{exams}[theo]{Examples}

\numberwithin{equation}{section}

\newcommand{\D }{\mathcal{D}}
\newcommand{\U }{\mathcal{U}}
\newcommand{\A }{\mathcal{A}}
\newcommand{\B }{\mathcal{B}}

\newcommand{\mA}{{\mathbb A}}
\newcommand{\mK}{{\mathbb K}}

\newcommand{\CC}{{\mathbb C}}

\newcommand{\Range}{{\rm Range }}

\newtheorem{theorem}{Theorem}[section]
\newtheorem{proposition}[theorem]{Proposition}
\newtheorem{corollary}[theorem]{Corollary}
\newtheorem{lemma}[theorem]{Lemma}
\newtheorem{definition}[theorem]{Definition}

\newtheorem{example}[theorem]{Example}
\newtheorem{remark}[theorem]{Remark}

\newcommand\cA{{\cal  A}}
\newcommand\cB{{\cal  B}}
\newcommand\cD{{\cal  D}}
\newcommand\cU{{\cal  U}}
\newcommand\cH{{\cal  H}}

\newcommand\cK{{\cal  K}}

\newcommand{\tr}{\,\mbox{\rm tr}}

\newcommand{\lb}{\label}

\newcommand{\bi}{\bibitem}

\newcommand{\beq}{\begin{equation}}
\newcommand{\eeq}{\end{equation}}

\title{
$2$-modified characteristic Fredholm determinants, Hill's method,
and the periodic Evans function of Gardner
}

\author{\sc \small
Kevin Zumbrun\thanks{Indiana University, Bloomington, IN 47405;
kzumbrun@indiana.edu: Research of K.Z. was partially supported
under NSF grant no. DMS-0300487.
 }}
\begin{document}

\maketitle


\begin{center}
\end{center}


\begin{abstract}
Using the relation established by Johnson--Zumbrun
between Hill's method of aproximating spectra of periodic-coefficient 
ordinary differential operators and a generalized periodic Evans function
given by the
$2$-modified characteristic Fredholm determinant of an associated
Birman--Schwinger system, together with a Volterra integral computation
introduced by Gesztesy--Makarov,
we give an explicit connection between the
generalized Birman--Schwinger-type periodic Evans function and
the standard Jost function-type
periodic Evans function defined by Gardner
in terms of the fundamental solution of the eigenvalue equation written
as a first-order system.
This extends to a large family of operators the results of 
Gesztesy--Makarov for scalar 
Schr\"odinger operators and 
of Gardner for vector-valued second-order elliptic operators,
in particular recovering by independent argument
the fundamental result of Gardner
that the zeros of the Evans function agree in location and (algebraic) 
multiplicity with the periodic eigenvalues of the associated operator.
\end{abstract}

\section{Introduction}

The purpose of this note, generalizing results of \cite{GM04,GLMZ05,GLM07}
for asymptotically constant-coefficient operators
and of Gesztesy--Makarov \cite{GM04} for
of periodic-coefficient Schr\"odinger operators, 
is to give an explicit connection between
two objects related to the spectra of periodic-coefficient 
operators on the line, namely, a Birman--Schwinger-type {\it characteristic
Fredholm determinant} $D(\lambda)$ introduced in \cite{JZ,BJZ} and the
Jost function-type {\it periodic Evans function} $E(\lambda)$
of Gardner \cite{G}.

Both of these objects are analytic functions whose zeros have
been shown in various contexts to agree in location and
multiplicity with the eigenvalues of the associated operator.
Thus, on the mutual domain for which these properties have been
established, it follows that they must agree up to a nonvanishing
analytic factor.
Here, we determine explicitly this nonvanishing factor, thus
illuminating the relation between the two functions while extending
the key property of agreement with eigenvalues to the union 
of the domains on which it has been established for each
function separately.
Knowledge of this scaling factor is
useful also for numerical approximation of spectra as
in \cite{BJZ,BJNRZ1,BJNRZ2}, allowing convenient comparison
of different methods.

Our method of proof proceeds by comparison with
an interpolating first-order Birman--Schwinger determinant,
using a simple version of the Volterra integral computation of \cite{GM04}
together with the relation established in \cite{JZ} between
Birman--Schwinger determinants and Hill's method, a spectral Galerkin algorithm
for finite approximation of spectra.
We carry out the analysis here for the simplest interesting
case of a second-order elliptic operator; however, our arguments
extend readily
to the full class of operators discussed in \cite{JZ}.

Consider the eigenvalue problem for a general second-order 
periodic-coefficient ordinary differential operator
$ L= B_0(x)^{-1}\big(\partial_x^2 + \partial_x A_1(x) + A_0(x)\big), $
written in form
\be\label{e:L}
\big(\partial_x^2 + \partial_x A_1(x) + A_0(x) -\lambda B_0(x)\big)U=0,
\ee
where $L$ is 
defined on complex vector-valued functions $U\in L^2[0,X]$ 
with periodic boundary conditions,
$A_j, B_j\in L^2$ are matrix-valued and periodic on $x\in [0,X]$,
and $B_0$ is positive or negative definite in the sense that its
symmetric part $\Re B_0:=(1/2)(B_0+B_0^*)$ is positive or negative
definite.\footnote{The assumption of symmetry of $B_0$ made in \cite{JZ} is
not necessary for the arguments there.}
(Note that, by Lyapunov's Lemma, 
definiteness can always be achieved by a variable--coefficient
change of coordinates,
provided that the eigenvalues of $B_0$ are of purely positive or purely
negative real part.)
Taking the Fourier transform, we may express \eqref{e:L}, equivalently,
as an infinite-dimensional matrix system
\be\label{infsys}
 (\D^2 +\D\A + \B)\U=0, 
\ee 
$ \B_{jk}= \hat A_{0,j-k} -\lambda \hat B_{0,j-k}$,
$\A_{jk}=\hat A_{1,j-k}$,
$\D_{jk}=(2\pi/X)\delta_j^k  ij$,
$\U_j= \hat U(j)$,
where $\hat f$ denotes the discrete Fourier transform of $f$.
(Here and elsewhere $i=\sqrt{-1}$.)
We shall alternate between these two representations as is convenient.

\section {The Birman--Schwinger Evans function}\label{s:review}
By a Birman--Schwinger-type procedure,
applying $(\partial_x^2-1)^{-1}$ to \eqref{e:L} on the left, 
we obtain an equivalent problem
\be\label{equiv}
 (I+K( \lambda))U=0,
\ee
where $K=K_1+K_0$, with 
$ K_1=\partial_x (\partial_x^2-1)^{-1} A_1,$
$K_0= (\partial_x^2-1)^{-1} (A_0+1 -\lambda B_0).$

\bl[\cite{JZ}]
For $A_j,B_j\in L^2$,
the operator $K$ is Hilbert-Schmidt.
\el

\begin{proof} (\cite{JZ})
Expressing $K$ in matrix form $\cK$ with respect to
the infinite-dimensional Fourier basis $\{e^{2\pi ijx/X}/\sqrt{X}\}$, 
we find that
$
\cK_{1,jk}= \frac{(2\pi i/X)j}{1+(2\pi/X)^2j^2} \hat A_{1,j-k},
$
and similarly for $\cK_2$,
where $\hat A_j$ denotes Fourier transform
and $i:=\sqrt{-1}$.
Taking without loss of generality $X=2\pi$,
we find that
\be\label{calc}
\|K_1\|_{\cB_2}=\|\cK_1\|_{\cB_2}=
\sum_j 
\frac{j^2}{(1+j^2)^2} \sum_k |A_1(j-k)|^2
=
\frac{j^2}{(1+j^2)^2} \|A_1\|_{L^2(x)}<+\infty,
\ee
and similarly for $\cK_2$.
\end{proof}  

\begin{definition}\label{BSEvans}
\textup{
We define the {\it Birman--Schwinger Evans function} as
\be\label{e:evans2}
D(\lambda):={\det}_2(I-K( \lambda)),
\ee
where $\det_2$ denotes the $2$-modified Fredholm determinant,
defined for Hilbert--Schmidt perturbations of the identity;
see Appendix \ref{s:fredholm} for a review of the relevant theory.
}
\end{definition}

\bpr[\cite{JZ,BJZ}]\label{anlem}
For $A_j,B_j \in L^2$, 
$D$ is analytic in $\lambda$; moreover, its zeros correspond in
location and multiplicity 
with the eigenvalues of $L:=B_0^{-1}(\partial_x^2+\partial_x A_1+ A_0)$. 
\epr

\begin{proof}
Analyticity follows by analyticity of the finite-dimensional Galerkin
approximations by which the Fredholm determinant is defined,
plus uniform convergence of the Galerkin approximations, a consequence
of \eqref{compare}.
Correspondence in location is immediate from the fact that
$\det_2(I-K)=0$ if and only if $(I-K)$ has a kernel, while
correspondence in multiplicity may be deduced by consideration
of a special sequence of Galerkin approximations on successive eigenspaces
of $L$, for which the Galerkin approximant can be seen to be a nonvanishing
multiple of the characteristic polynomial for the restriction of $L$
to these finite-dimensional invariant subspaces.
See \cite{JZ} for further details.
\end{proof}

\section{Hill's method}\label{hills}

Hill's method consists of truncating \eqref{infsys} at wave number
$J$, i.e., considering the $(2j+1)$-dimensional minor $|j|\le J$,
and solving the resulting finite-dimensional system to obtain
approximate eigenvalues for $L$: that is, the eigenvalues of
the $(2J+1)\times (2J+1)$ matrix
\be\label{LJ}
 L_J=(\B_{0,J})^{-1}(\D_J^2 +\D_J\A_J + \B_J) ,
\ee
$\B_{0,J,jk}=\hat B_{0,j-k}$,
where subscripts $J$ indicate truncation at wave number $J$, or, equivalently,
restriction to the $J$th centered minor as described above.

Left-multiplying by $(\cD_J^2-1)^{-1}$, similarly as above,
we may rewrite the truncated eigenvalue equation
 $(L_J-\lambda)U=(\D_J^2 +\D_J\A_J + \B_J)\cU=0 $ equivalently as
$ (I+ \cK_J)\cU=0, $
where $\cK_J=\cK_{1,J}+\cK_{2,J}$ is the truncation of the
Fourier representation $\cK=\cK_1+\cK_2$ of operator $K$, that is,
 $\cK_1=\D_J(\D_J^2-I)^{-1}\A_{1,J}$,
 $\cK_2=(\D_J^2-I)^{-1}(\A_{0,J}+1-\lambda \cB_{0,J})$.
We define the
{\it truncated Birman--Schwinger Evans function}, accordingly, as 
\be\label{truncevans}
D_J(\lambda):={\det}_2(I-\cK_J)=
\det (\cD_J^2-1) \det \cB_{0,J} e^{\tr \cK_J(\lambda)} \det (L_J-\lambda).
\ee

\bpr[\cite{JZ}]\label{trunccorr}
Each $D_J$ is analytic in $\lambda$; moreover,
the zeros of $D_{J}$ correspond
in location and multiplicity with those of $L_{J}$.
\epr

\begin{proof}
Immediate,
from \eqref{truncevans}, the assumed uniform positive definiteness
of $B_0$, hence also of $\cB_{0}$ and $\cB_{0,J}$, 
and properties of the characteristic polynomial.
\end{proof}


\bpr[\cite{JZ}] \label{convthm}
For $A_j, B_j\in L^2$, $D_{J}(\lambda)\to D(\lambda)$ 
as $J\to \infty$, uniformly on $|\lambda|\le R$.
\epr

\begin{proof}
By diagonality of $\cD$, $\cD_J\cA_J=
(\cD\cA)_J$, 
hence $D_{J}$ is a sequence of Galerkin 
approximations of the Fredholm determinant $D=\det_2(I-\cK)$ 
on a complete set of successively larger subapaces.
By definition of the Fredholm determinant, therefore,
$D_{J}\to D$ as $J\to \infty$.
\end{proof}

\bc 
[Convergence of Hill's method \cite{JZ}] \label{Lconvthm}
For $A_j, B_j\in L^2$, 
 the eigenvalues of $L_{J}$
approach the eigenvalues of $L$ in location and multiplicity
as $J\to \infty$, uniformly on $|\lambda|\le R$.
\ec

\begin{proof}
Immediate, from Lemma \ref{anlem}, Propositions \ref{trunccorr} 
and \ref{convthm}, and properties of uniformly convergent analytic functions.
\end{proof}

\section{The Evans function of Gardner}

Rewriting \eqref{e:L} as a first-order system
\be\label{e:1L}
\partial_x W=A(\lambda)W,
\quad
W:=\bp U\\ U_x\ep,
\quad
\mA(\lambda):=\bp 0 & I\\
(\lambda B_0-A_0 -\partial_x A_1) & -A_1
\ep,
\ee
we may alternatively define an Evans function by 
shooting, as done by Gardner.

\begin{definition}\label{def:GEvans}
\textup{
We define the {\it Jost-type Evans function} of Gardner, following \cite{G}, as
\be\label{e:GEvans}
E(\lambda):=
\det (\Psi(X)-I),
\ee
where $\Psi$ denotes the fundamental solution of \eqref{e:1L},
satisfying 
$\Psi(0)=I$, and $X$ is the period.
}
\end{definition}

\bpr[\cite{G}]\label{gardlem}
For $A_j,B_j \in C^1$, 
$D$ is analytic in $\lambda$; moreover, its zeros correspond in
location and multiplicity 
with the (periodic) eigenvalues of $L:=B_0^{-1}(\partial_x^2+\partial_x A_1+ A_0)$. 
\epr

\begin{proof}
Analyticity and agreement in location follow immediately 
by analytic dependence of solutions of ODE and the fact that $\lambda$ 
is an eigenvalue precisely if there exists a solution 
$W(x)=\Psi(x)W_0$ of \eqref{e:1L} for which $W(X)=\Psi(X)W_0= W_0$.
Agreement in multiplicity may be shown by a more detailed argument
based on choice of a Jordan basis as in \cite{G}.
\end{proof}

\section{Connections}

Define now the intermediate
{\it first-order Birman--Schwinger Evans function}
\be\label{int}
F(\lambda):={\det}_2
\big( (\partial_x-1)^{-1}(\partial_x-\mA)\big)
={\det}_2(I-\mK), 
\ee
where $\mK:=(\partial_x-1)^{-1}(\mA - 1)$.

\bt\label{conn1}
For $A_j,B_j \in L^2$, 
\be\label{FE}
F(\lambda)=
\gamma (e^X-1)^{-2n} E(\lambda),
\ee
where 
\be\label{e:gamma}
\gamma:= e^{ \frac{e^X}{1-e^X}
\int_0^X(\tr(\mA)-2n)(y)dy}=e^{\frac{e^X}{e^X-1}(\tr(A_{1,ave})+2n)X},
\ee
and $A_{1,ave}$ denotes the mean over one period of $A_1$.
\et

\br\label{Xinf}
\textup{
From \eqref{FE}-\eqref{e:gamma}, we obtain
$
F(\lambda)\sim 
E(\lambda) e^{ -\int_0^X \tr(\mA)(y)dy}
=
\det(I- \Psi(X)^{-1})
$
as $X\to \infty$, so that $F$ is asymptotic to a ``backward''
version of Gardner's Evans function.
}
\er

\begin{proof}
By a straightforward computation taking into account periodicity 
together with the jump condition $G(y^+,y)-G(y_-,y)=1$ at $x=y$, 
we find that the Green kernel associated with $(\partial_x-1)^{-1}$ is
semiseparable, of form
\be\label{G}
G(x,y)=\begin{cases}
\frac{e^{x-y}}{1-e^X},&x>y,\\
\frac{ e^{X+x-y}}{1-e^X},&x<y.\\
\end{cases}
\ee

Following \cite{GM04}, we may thus decompose $G$ as the sum
$G(x,y)=f(x)g(y)+J(x,y)$ of a separable kernel 
\be
f(x)g(y):\quad
f(x)= e^x,
\quad
g(y)= \frac{e^X e^{-y}}{1-e^X},
\ee
and a Volterra kernel
\be\label{H}
J(x,y)=
\begin{cases}
e^{x-y},&x>y,\\
0,&x<y.\\
\end{cases}
\ee
Likewise, the kernel $\kappa(x,y)$ of $\mK$ decomposes as
\be\label{decomp}
\kappa(x,y)=f(x)g(y)(\mA(y)-1) +H(x,y), 
\ee
where $H(x,y=J(x,y)(\mA(y)-1)$ is again a Volterra kernel,
vanishing for $x<y$.

Using the Fredholm determinant facts \eqref{mult}, \eqref{commute}
compiled in Appendix \ref{s:fredholm},
together with the fact that ${\det}_2(I-H)=1$ for a Volterra operator
$H$,\footnote{
Heuristically a strictly lower triangular perturbation of the identity;
$H$ has no nonzero eigenvalues.}
we find (here, with slight abuse of notation, denoting operators 
by their associated kernels) that
\ba\label{reduct1}
{\det}_2 \big(I- &(\partial_x-1)^{-1}(\mA-1)\big)\\
& = {\det}_2 \big(I- H -fg(\mA-1))\\
&= {\det}_2 \big(I- H\big) {\det}_2 \big(I- (I- H)^{-1}fg(\mA-1)\big)
e^{-\tr\big(H(I-H)^{-1}fg(\mA-1) \big)}\\
&= {\det}_2 \big(I- (I- H)^{-1}fg(\mA-1)\big)
e^{\tr\big(fg(\mA-1) \big)}
e^{-\tr\big((I-H)^{-1}fg(\mA-1) \big)}\\
&= \gamma {\det}_{\CC^{2n}} \big(I- \langle g(\mA-1),(I- H)^{-1}f \rangle \big),
\ea
where $\langle \cdot, \cdot\rangle$ denotes $L^2[0,X]$ inner
product and 
$\gamma:= e^{\tr\big(\langle g(\mA-1), f\rangle \big)}=
e^{\frac{e^X}{e^X-1}(\tr(A_{1,ave})+2n)X}$.

Computing 
$$
(I-H)\Psi(x)= \Psi(x) - \int_0^x e^{x-y}(\mA-1)\Psi(y)dy,
$$
and recalling that $(\partial_x -1)\Psi=(\mA-1)\Psi$, where
$\Psi$ denotes the fundamental solution of $(\partial_x -\mA)\Psi=0$,
$\Psi(0)=I$,
so that, by Duhamel's principle/variation of constants,
$\Psi(x)=e^x + \int_0^x e^{x-y}(\mA-1)\Psi(y)dy$,
we find that
$(I-H)\Psi= e^x =f$, or $(I-H)^{-1}f=\Psi$.
Substituting into \eqref{reduct1}, we obtain,
noting that $\Theta:=e^{-y}\Psi(y)$ satisfies
$\partial_y \Theta= (\mA-1)\Theta$,
\ba\label{finalred}
{\det}_2 \big(I- (\partial_x-1)^{-1}(\mA-1)\big)&=
 \gamma{\det}_{\CC^{2n}} \Big(I- \Big(\frac{e^X}{1-e^X}\Big)
\int_0^X(\mA-1) e^{-y}\Psi(y)dy \Big)\\
&= \gamma {\det}_{\CC^{2n}} 
\Big(I- \Big(\frac{e^X}{1-e^X}\Big) e^{-y}\Psi(y)|_0^X \Big)\\
&=\gamma {\det}_{\CC^{2n}} 
\Big(I- \Big(\frac{e^X}{1-e^X}\Big) (e^{-X}\Psi(X)-I) \Big)\\
&= \gamma (e^X-1)^{-2n} {\det}_{\CC^{2n}} (\Psi(X)-I),
\ea
yielding the result.
\end{proof}

Now, define the
{\it truncated first-order Birman--Schwinger Evans function} as 
\be\label{1truncevans}
F_J(\lambda):={\det}_2(I-\hat \cK_J),
\ee
where $\hat \cK_J$ denotes the $J$th Fourier truncation of the
Fourier transform $\hat \cK$ of $\mK$.
\bl\label{oneconv}
For $A_j, B_j\in L^2$, $F_{J}(\lambda)\to F(\lambda)$ 
as $J\to \infty$, uniformly on $|\lambda|\le R$.
\el

\begin{proof}
As in the argument of Proposition \ref{convthm}, 
this follows essentially by the definition of the Fredholm determinant.
\end{proof}

\bl\label{sim}
Similarly as in \eqref{truncevans}, we have
\be\label{eval1truncevans}
F_J(\lambda)
= \det(\cD_J-1)^2
\det \cB_{0,J} e^{\tr \hat \cK_J(\lambda)} \det (L_J-\lambda).
\ee
\el

\begin{proof}
Using $(I-\mK)=(\partial_x-1)^{-1}(\partial_x - \mA )$, 
$\mA(\lambda)=\bp 0 & I\\ (\lambda B_0-A_0 -\partial_x A_1) & -A_1 \ep $
to obtain
$$
(I-\hat \cK_J)=\bp \cD_J-I_J & 0\\
0& \cD_J-I_J \ep^{-1}
\bp
\cD_J & -I_J\\
\cB_J +\widehat {(\partial_x A_1)}_J &\cD_J +\cA_J
\ep,
$$
and, using the block determinant formula
$\det \bp a & -I\\b&c\ep=
\det(ca+b)$,
$$
\begin{aligned}
\det \bp
\cD_J & -I_J\\
\cB_J +\widehat {(\partial_x A_1)}_J &\cD_J +\cA_J
\ep
&=
\det\big((\cD_J+\cA)\cD_J +\widehat {(\partial_x A_1)}_J+ \cB_J \big)\\
&=
\det\big(\cD_J^2+ \cD_J \cA_J + \cB_J \big)\\
&=
\det \cB_{0,J} \det (L_J-\lambda).
\end{aligned}
$$
\end{proof}

\bc\label{conn2}
For $A_j,B_j \in L^2$, $ F(\lambda)= \epsilon 
e^{\hat \delta- \delta} D(\lambda) , $
where, for $X=2\pi$,
\ba\label{deltagamma}
\hat \delta&:= \tr (\hat A_{1,0}+2I)\Big(1+ \sum_{1\le |j|} 
\frac{1}{j^2+ij}\Big),\\
\delta&:= -\tr(\hat A_{1,0}+I -\lambda \hat B_{0,0})\sum_{j}\frac{1}{j^2+1},\\
\epsilon&:=\lim_{J\to \infty}\prod_{|j|\le J} \Big(1+\frac{2}{ij+1}\Big).
\ea
\ec

\begin{proof}
We may readily calculate that
\ba\label{trK}
\tr \cK&= \tr(\hat A_{1,0})\sum_{|j|\le J}\frac{j}{-j^2-1}
+
\tr(\hat A_{1,0}+I -\lambda \hat B_{0,0})\sum_{|j|\le J}\frac{1}{-j^2-1}
\\
&=
-\tr(\hat A_{1,0}+I -\lambda \hat B_{0,0})\sum_{|j|\le J}\frac{1}{j^2+1}
\ea
and
\ba\label{trhatK}
\tr \hat \cK&=
\bp -(\cD_J-I_J)^{-1} & I\\ 
*
&
(\cD_J-I_J)^{-1}(-A_1-I) \ep \\
& = -\tr (\hat A_{1,0}+2I)\sum_{|j|\le J} \frac{1}{ij-1}\\
& = -\tr (\hat A_{1,0}+2I)\Big(-1+ \sum_{1\le |j|\le J} 
\Big( \frac{1}{ij-1}-\frac{1}{ij}\Big)\Big)\\
& = \tr (\hat A_{1,0}+2I)\Big(1+ \sum_{1\le |j|\le J} 
\frac{1}{j^2+ij}\Big),
\\
\ea
both uniformly convergent by absolute convergence of $\sum 1/j^2$.
Similarly, 
$$
\prod_{|j|\le J}\frac{(ij-1)^2}{-j^2-1}=
\prod_{|j|\le J} \Big(1+\frac{2}{ij+1}\Big)
\sim e^{\sum_{|j|\le J} \frac{2}{ij+1} + O(j^{-2})}
\sim e^{\sum_{|j|\le J} O(j^{-2})}
$$
may be seen to be uniformly convergent.
Comparing 
\eqref{truncevans} and \eqref{eval1truncevans}
and taking the limit as $J\to \infty$, we obtain the result.
\end{proof}


Collecting information, we have our final result, relating
$D$, $E$, and $F$.

\bc\label{finalconn}
For $A_j,B_j \in L^2$, and $X=2\pi$,
\be\label{finalrel}
D(\lambda)=
\frac{e^{\delta-\hat \delta}}{\epsilon}
F(\lambda) =
\frac{e^{\delta-\hat \delta}}{\epsilon} \gamma (e^X-1)^{-2n} E(\lambda),
\ee
where $\gamma$ is as in \eqref{e:gamma}
and $\delta$, $\hat \delta$, and $\epsilon$
are as in \eqref{deltagamma}.
\ec




\br\label{gardrmk}
\textup{
Since $\frac{e^{\delta-\hat \delta}}{\epsilon} \gamma (e^X-1)^{-2n}$
is analytic in $\lambda$ and nonvanishing, 
Corollary \ref{finalconn} together with
Proposition \ref{anlem} gives an alternative 
proof of Proposition \ref{gardlem}.
}
\er

\appendix
\section {$2$-modified Fredholm determinants}\label{s:fredholm}
In this appendix, we recall for completeness the basic properties
of $2$-modified Fredholm determinants of 
Hilbert--Schmidt perturbations of the identity
\cite{GGK97,GGK00,GK69,Si77,Si05}.
For a Hilbert space $\cH$, the Hilbert--Schmidt class
$\cB_2(\cH)$ is defined as the set of linear operators on $\cH$
for which 
$$
\|A\|_{\cB_2(\cH)}
:=\sum_{j,k} |\langle Ae_j,e_k\rangle|^2=\tr ( A^*A)<+\infty,
$$
where $\{e_j\}$ is any orthonormal basis.
Evidently, $\|\cdot\|_{\cB_2(\cH)}$ is independent of the basis.

For a finite-rank operator $A$, the $2$-modified Fredholm determinant 
is defined as 
\ba \lb{2.34} 
 {\det}_{2,\cH} (I_{\cH}-A):= {\det}_{\cH}((I_{\cH}-A)e^{A})
={\det}_{\cH}(I_{\cH}-A) \, e^{\tr_{\cH}(A)},   
\ea
where $\det_{\cH}$ denotes the usual determinant restricted to $\Range(A)$.
There holds
	\be \lb{condition} 
e^{-C\|A\|_{\cB_2(\cH)}^2}\le
 |{\det}_{2,\cH}(I_{\cH}-A)| \leq e^{C\|A\|_{\cB_2(\cH)}^2}
\ee
and
\be\label{compare}
 |{\det}_{2,\cH}(I_{\cH}-A) - {\det}_{2,\cH}(I_{\cH}-B)| \leq \|A-B\|_{\cB_2(\cH)}
e^{C[\|A\|_{\cB_2(\cH)}+\|B\|_{\cB_2(\cH)}+1]^2},
\ee
where $C>0$ is a constant independent of the dimension of the space.

For $A\in \cB_2(\cH)$, the
$2$-modified Fredholm determinant is defined as the limit
\be\label{limdef}
 {\det}_{2,\cH}(I_{\cH}-A):= \lim_{J\to \infty} {\det}_{2,\cH_J}(I_{\cH_J}-A_J),
\ee
where $\cH_J$ is any increasing sequence of finite-dimensional subspaces
filling up $\cH$, and $A_J$ denotes the Galerkin approximation
$P_{\cH_J}A|_{\cH_J}$, where $P_J$ is the orthogonal projection
onto $\cH_J$.
Equivalently,
 ${\det}_{2,\cH}(I_{\cH}-A):= \Pi_j (1-\alpha_j)e^{\alpha_j} $,
where $\alpha_j$ are the eigenvalues of $A$. 
For $A\in \cB_2$, $(I-A)$ is invertible if and only if
$\det_2(I-A)\ne 0$.

Properties \eqref{2.34}--\eqref{condition} are evidently 
inherited by continuity.  
Likewise, from the corresponding properties
of finite-dimensional determinants, we obtain in the limit
\be\label{mult}
 {\det}_{2,\cH}((I_{\cH}-A)(I_{\cH}-B))={\det}_{2,\cH}(I_{\cH}-A)
{\det}_{2,\cH}(I_{\cH}-B) \, e^{-\tr_{\cH}(AB)}
\ee
for all $A,B\in \cB_2$ (here we are implicitly
using the fact that $AB$ is
in trace class, with $\|AB\|_{\cB_1}\le \|A\|_{\cB_2} \|B\|_{\cB_2}$,
so that $\tr (AB)$ is well-defined)
 and
\be\label{commute}
{\det}_{2,\cH'}(I_{\cH'}-AB)={\det}_{2,\cH}(I_\cH-BA)
\ee
for all $A\in\cB(\cH,\cH')$, $B\in\cB(\cH',\cH)$
such that $BA\in \cB_2(\cH)$, $AB\in \cB_2(\cH')$.

\end{document}